\documentclass[10pt]{article}

\usepackage[affil-it]{authblk}
      \usepackage{latexsym}
         \usepackage[reqno, namelimits, sumlimits]{amsmath}
         \usepackage{amssymb, amsfonts}
         \usepackage{amsthm}

 \newtheorem{theorem}{Theorem}[section]
 
 \newtheorem{definition}[theorem]{Definition}
 
 \newtheorem{lemma}[theorem]{Lemma}
 
 \newtheorem{remark}[theorem]{Remark}

 \newtheorem{pro}[theorem]{Proposition}
 
\title{  Local boundary regularity for the Navier-Stokes equations in nonendpoint borderline Lorentz spaces}
\author{T Barker
  \thanks{Email address: \texttt{tobias.barker@seh.ox.ac.uk}; }}
\affil{OxPDE, Mathematical Institute, University of Oxford, Oxford,UK}
\date{ \today}

\begin{document}
\maketitle
\begin{abstract}
\end{abstract}
\setcounter{equation}{0}
We prove local regularity up to flat part of boundary, for certain classes of distributional solutions  that are $L_{\infty}L^{3,q}$ with $q$ finite. 
\setcounter{equation}{0}
\section{Introduction}
  In this paper we are going to prove local regularity, up to flat part of the boundary boundary, for  certain classes of weak solutions to the three dimensional Navier Stokes equations. The main assumption is the velocity field belongs to $L_{\infty}L^{3,q}$ with $q$ finite.\\
In the local theory, two cases are distinguished: interior and boundary. The first result, regarding local interior regularity criteria for the Navier-Stokes in terms of the velocity, was established by Serrin in \cite{Serrin}. Later the following generalisation was proven by Struwe in \cite{Struwe}.
\begin{theorem}\label{struwe}
Let $v$ be a divergent free vector field defined in the unit parabolic cylinder $Q=B\times ]-1,0[$, where $B$ is the unit ball in $\mathbb{R}^3$ centered at the origin. Suppose that $v$ satisfies the three conditions
\begin{equation}\label{dissipation/energy}
v\in L_{2,\infty}(Q)\cap W^{1,0}_2(Q),
\end{equation}
\begin{equation}\label{weaksolution}
\int\limits_{Q} v.\partial_{t}\phi- v\otimes v:\nabla \phi+\nabla v:\nabla \phi dz=0
\end{equation}
for all smooth  solenodial functions $\phi$ that are compactly supported in $Q$ (denoted $C^{\infty}_{0,0}(Q)$,
and
\begin{equation}\label{LPS}
v\in L_{s,l}(Q),\, 3/s+2/l=1,\, s>3.
\end{equation}
\end{theorem}
Before further comment, we give the required notation. 
If $X$ is a Banach space with norm $\|\cdot\|_{X}$, then $L_{s}(a,b;X)$, $a<b$, means the usual Banach space of strongly measurable $X$-valued functions $f(t)$ on $]a,b[$ such that the norm
$$\|f\|_{L_{s}(a,b;X)}:=\left(\int\limits_{a}^{b}\|f(t)\|_{X}^{s}dt\right)^{\frac{1}{s}}<+\infty$$
for $s\in [1,\infty[$, and with the usual modification if $s=\infty$. With this notation if $Q_{T}=\Omega\times ]0,T[$ then 
$$L_{s,l}(Q_{T}):= L_{l}(0,T; L_{s}(\Omega)).$$
We  define the following Sobolev space with the mixed norm:
$$ W^{1,0}_{m,n}(Q_{T})=\{ v\in L_{m,n}(Q_{T}): \|v\|_{L_{m,n}(Q_{T})}+$$$$+\|\nabla v\|_{L_{m,n}(Q_{T})}<\infty\}.$$
For $s=3$ the following local interior regularity result  was proven in \cite{ESS2003}. Namely the following:

\begin{theorem}\label{ESS}
Suppose that a pair of functions $(u,p)$ satisfies the Navier-Stokes equations in $Q({1})$ in the sense of distributions such that (for 
\begin{equation}\label{uenergyclass}
u\in L_{2,\infty}(Q)\cap L_{2}(-1,0;W^{1}_{2}(B)),
\end{equation}
and
\begin{equation}\label{pressurerestriction}
p\in L_{\frac 3 2}(Q).
\end{equation}
Suppose further that
\begin{equation}\label{bddL3norm}
u\in L_{3,\infty}(Q)
\end{equation}
for some $q\in ]3,\infty[$. Then the velocity function $u$ is Holder continuous on $\bar{Q}({1}/{2})$.

\end{theorem}
Here we define:
$$B(x_0,R)=\{x\in\mathbb{R}^3: |x-x_0|<R\},$$
$$B(\theta)=B(0,\theta),\,\,\,B=B(1),$$
$$ Q(z_0,R)=B(x_0,R)\times ]t_0-R^2,t_0[,\,\,\, z_{0}=(x_0,t_0),$$
$$Q(\theta)=Q(0,\theta).$$
Theorem \ref{ESS} is was proven in \cite{ESS2003} by ad absurdum. A suitable rescaling and limiting procedure is performed and gives a nontrivial solution to Navier-Stokes equations in $\mathbb{R}^3\times ]0,\infty[$. Then contradiction is then obtained by proving a Liouville type theorem involving backward uniqueness for a certain class of parabolic operators. Subsequently, a version of Theorem \ref{ESS} was proven up to the flat part of the boundary in \cite{S2005}. This was also done near the curved part of the boundary in \cite{MSh2006}. In the context of Lorentz spaces, the interior result is proven in \cite{WangZhang1}.
Namely the assumptions and statement proved  are the same as in Theorem \ref{ESS} except (\ref{bddL3norm}) is replaced by
\begin{equation}\label{bddLorentznorm}
u\in L_{\infty}(-1,0; L^{3,q}(B))
\end{equation}
with $q$ finite. Recently, a version is proven in \cite{Phuc} but with the additional restriction that 
\begin{equation}\label{Phucpres}
p\in L_{2}(-1,0; L_{1}(B)).
\end{equation}
We will first present a different proof to \cite{WangZhang1} of Theorem \ref{ESS} under the assumption that $u\in L_{\infty}(-1,0; L^{3,q}(B)).$ The contradiction argument has the same spirit as that used in \cite{ESS2003}, \cite{Phuc} and \cite{WangZhang1}.
The difference  is the development of new estimates of certain scale invariant quantities associated with the pressure and velocity, these may be of independent interest. Our reasons for first presenting the known interior result are two fold. Firstly the estimates of scale invariant quantities used there mostly carry (with some adjustment) to the case involving the flat part of the boundary. Secondly, the  interior result is a necessary prerequisite for proving the boundary case.\\\\
Now we can state our main goal of local regularity up to flat part of the boundary for nonendpoint Lorentz spaces.
\begin{theorem}\label{Barkerboundaryreg}
Let a pair of functions $v$ and $p$ have the following differentiability properties:
\begin{equation}\label{updifprop}
v\in L_{2,\infty}(Q^+(2))\cap W_{2}^{1,0}(Q^+(2)),\,\,\,\, p\in L_{\frac{3}{2}}(Q^+(2)).
\end{equation}

Suppose that $v$ and $p$ satisfies the Navier-Stokes equations
\begin{equation}\label{navierstokeseqnhalf}
\partial_{t}v+\rm{div}\,v\otimes v-\Delta v=-\nabla p,\,\,\rm{div}\,v=0
\end{equation}
 in $Q^{+}(2)$ in the sense of distributions along with the boundary condition
\begin{equation}\label{tracezero}
v(x,t)=0,\,\,\, x\in\Gamma(2)\,\, and\,\, -4< t< 0.
\end{equation}

Assume, in addition, that there exists  $3\leqslant q<\infty$ such that
\begin{equation}\label{Lorentzhalfbounded}
v\in L_{\infty}(-4,0; L^{3,q}(B^{+}(2)).
\end{equation}
Then $v$ is Holder continuous in the closure of the set
$Q^+({1}/{2}).$
\end{theorem}
 We explain the notation. 
Setting $x'=(x_1,x_2)\in\mathbb{R}^2$, we introduce the following definitions:

$$B^{+}(x_0,R)=\{x\in B(x_0,R):x=(x',x_3),\,\,x_3>x_{03}\},$$
$$B(\theta)=B(0,\theta),\,\,\,B=B(1)\,\,\,B^+(\theta)=B^+(0,\theta),\,\,\,B^{+}=B^+(1),$$
$$\Gamma(x_0,R)=\{x\in B(x_{0},R): x_{3}=x_{30}\},\,\,\, \Gamma(\theta)=\Gamma(0,\theta),\,\,\Gamma=\Gamma(1),$$

$$Q^+(z_0,R)=B^+(x_0,R)\times ]t_0-R^2,t_0[,$$

$$Q^+(\theta)=Q^+(0,\theta),\,\,\,Q^+=Q^+(1).$$ 
Before commenting further let us define the Lorentz spaces. 
For a measurable function $f:\Omega\rightarrow\mathbb{R}$ define:
\begin{equation}\label{defdist}
d_{f,\Omega}(\alpha):=|\{x\in \Omega : |f(x)|>\alpha\}|.
\end{equation}
Given a measurable subset $\Omega\subset\mathbb{R}^{n}$,  the Lorentz space $L^{p,q}(\Omega)$ (with $p\in ]0,\infty[$, $q\in ]0,\infty]$) is the set of all measurable functions $g$ on $\Omega$ such that the quasinorm $\|g\|_{L^{p,q}(\Omega)}$ is finite. Here:

\begin{equation}\label{Lorentznorm}
\|g\|_{L^{p,q}(\Omega)}:= \Big(p\int\limits_{0}^{\infty}\alpha^{q}d_{g,\Omega}(\alpha)^{\frac{q}{p}}\frac{d\alpha}{\alpha}\Big)^{\frac{1}{q}},
\end{equation}
\begin{equation}\label{Lorentznorminfty}
\|g\|_{L^{p,\infty}(\Omega)}:= \sup_{\alpha>0}\alpha d_{g,\Omega}(\alpha)^{\frac{1}{p}}.
\end{equation}\\
It is well known that for $q\in ]0,\infty[,\,q_{1}\in ]0,\infty]$ and $q_{2}\in ]0,\infty]$ with $q_{1}< q_{2}$ we have the  embedding
$ L^{p,q_1} \hookrightarrow  L^{p,q_2}$
 and the inclusion is known to be strict.
Roughly speaking, the second index of Lorentz spaces gives information regarding nature of logarithmic bumps. For example, using decreasing rearrangements, it can be verified that for any $1>\beta>0, q>3$ we have
\begin{equation}\label{lorentzlogsing}
|x|^{-1}|\log(|x|^{-1})|^{-\beta}\chi_{|x|< 1}(x)\in L^{3,q}(\mathbb{R}^3)\,\,\,\rm{if\, and\, only\, if}\,\,q>\frac{1}{\beta}.
\end{equation}
 In this way Theorem \ref{Barkerboundaryreg} gives a strengthening of the previous result obtained in \cite{S2005}.  It should be stressed that, at the time of writing, the question of interior regularity  is open for the critical norm $L^{3,\infty}(B)$ that contains $|x|^{-1}$. We mention that  interior regularity results, that have smallness condition on $L_{\infty}(L^{3,\infty})$ norm, have been obtained in  \cite{kimkozono}, \cite{kozono} and \cite{Tsai}, for example. 
\\Now we can explicitly describe the challenges presented by the boundary and  to motivate the method used. The proof also proceeds by contradiction and uses the proof of a Liouville type theorem via backward uniqueness. The major differences lie in the treatment of the pressure. Unlike the interior case, if we scale and blow up the Navier-Stokes equations at singular boundary points we do not know if we can obtain a limiting pressure for the boundary case that lies in the space $L_{\infty}(-\infty,0; L^{\frac{3}{2},\frac{q}{2}}(\mathbb{R}^{3}_{+}))$. Unfortunately, we cannot even show that there is a reasonable global norm of the limiting pressure
which is finite.
 In our investigation, we were only able to demonstrate that the limiting pressure has the same integrability, on compact space-time subsets of $\mathbb{R}^{3}_{+}\times ]-\infty,0[$, as that for the original pressure. This creates major difficulties with the epsilon regularity used in the the interior regularity result in Lorentz spaces, presented in \cite{Phuc}. This criteria needs the limiting pressure to have more local integrability in time than that assumed for the original pressure in Theorem \ref{Barkerboundaryreg}. Our investigation is in the same spirit as that of \cite{S2005}.
 The main differences are that we have to strengthen Lemmas on the decay of the pressure (Proposition 2.5-2.6 of that paper), together with the development of a convenient epsilon regularity criteria for interior regularity of suitable weak solutions. Both of these points use the scale invariant estimates used in our version of the proof of the corresponding interior result.  The epsilon regularity criteria may be of independent interest, as it provides a strengthening to the statement given in \cite{Ladyser02} and \cite{Lin}.


\setcounter{equation}{0}

\section{Local interior regularity in nonendpoint borderline Lorentz spaces}

Here is the explicit statement, which we provide a proof of. An alternative proof can be found in \cite{WangZhang1}.
\begin{theorem}\label{Interiorreg}
Suppose that a pair of functions $(u,p)$ satisfies the Navier-Stokes equations in $Q({1})$ in the sense of distributions such that (for 
\begin{equation}\label{uenergyclass1}
u\in L_{2,\infty}(Q)\cap L_{2}(-1,0;W^{1}_{2}(B)),
\end{equation}
and
\begin{equation}\label{pressurerestriction1}
p\in L_{\frac{3}{2}}(Q).
\end{equation}
Suppose further that
\begin{equation}\label{bddLorentznorm1}
u\in L_{\infty}(-1,0; L^{3,q}(B))
\end{equation}
for some $q\in ]3,\infty[$. Then the velocity function $u$ is Holder continuous on $\bar{Q}({1}/{2})$.

\end{theorem}
We briefly recap the  definition of a suitable weak solution to the Navier-Stokes equations given by Lin in \cite{Lin}. Here it is.
\begin{definition}\label{suitweakinter}
Let $\omega$ be an open set in $\mathbb{R}^{3}$. We say that a pair $v$ and $p$ are a suitable weak solution to the Navier-Stokes equations on the set $\omega\times ]-T_{1},T[$ if they satisfy the Navier Stokes equations in the sense of distributions. Moreover they are required to satisfy the following conditions.
\begin{equation}\label{vdifsuitinter}
v\in L_{2,\infty}(\omega\times ]T_1,T[)\cap L_{2}(-T_{1},T; W^{1}_{2}(\omega)),
\end{equation}
\begin{equation}\label{pdifsuitinter}
p\in L_{\frac{3}{2}}(\omega\times ]-T_{1},T[).
\end{equation}

For a.a $t\in ]-T_1,T[$ and for all non negative cut-off functions $\phi\in C^{\infty}_{0}(\mathbb{R}^4)$ vanishing in a neighbourhood of the parabolic boundary
$$\partial^{'}Q= \omega\times \{t=-T_{1}\}\cup\partial \omega\times [-T_{1},T],$$
 $v$ and $p$ satisfy the local energy inequality
\begin{equation}\label{localenergyinequalityinter}
\int\limits_{\omega}\phi(x,t)|v(x,t)|^2dx+2\int\limits_{\omega\times ]-T_{1},t[}\phi |\nabla v|^2 dxdt^{'}\leqslant$$$$\leqslant
\int\limits_{\omega\times ]-T_{1},t[}[|v|^2(\partial_{t}\phi+\Delta\phi)+v.\nabla\phi(|v|^2+2p)] dxdt^{'}.
\end{equation}
\end{definition}
Let us proceed with a Lemma. The analogous Lemma (Lemma 4.1) was stated and proven in \cite{Phuc}. As the proof of this statement is essentially unchanged we omit it.
\begin{lemma}\label{usuitinter}
Suppose that the pair of functions $(u,p)$ satisfy the hypothesis of Theorem \ref{Interiorreg}. Then $(u,p)$ forms a suitable weak solution to Navier-Stokes equations in $Q({5}/{6})$ with a generalized energy inequality and furthermore $u\in L_{4}(Q)$.
Moreover the inequality
\begin{equation}\label{uLinfinitylorentz}
\|u(\cdot,t)\|_{L^{3,q}(B(3/4))}\leqslant \|u\|_{L_{\infty}(-(3/4)^{2},0;L^{3,q}(B(3/4))},
\end{equation}
holds for all $t\in [-(3/4)^2,0]$, and the function 
$$t\rightarrow\int\limits_{B({3}/{4})}u(x,t)w(x)dx$$
is continuous on $[-({3}/{4})^{2},0]$ for any $w\in L^{\frac{3}{2},\frac{q}{q-1}}(B({3}/{4})).$
\end{lemma}

Next, we recap the rescaling prodecure for the proof of Theorem \ref{Interiorreg} used by Escauriza, Seregin, Sverak in \cite{ESS2003}.
Suppose the conditions for $(u,p)$ of Theorem \ref{Interiorreg} hold.
Then by the previous Lemma $(u,p)$ form a suitable weak solution to the Navier-Stokes equations in $Q({5}/{6})$ and for $t\in [-({3}/{4})^2,0]$ we have
\begin{equation}\label{pointwiseLorentz}
u(\cdot,t)\in L^{3,q}(B({3}/{4})).
\end{equation}
The rescaling procedure arises from assuming Theorem \ref{Interiorreg} is false. Thus, $u$ has no representative that is Holder continuous on $\bar{Q}({1}/{2})$. This implies that there exists a singular point $z_{0}\in\bar{Q}({1}/{2})$ such that there is no parabolic neighbourhood of $O_{z_{0}}$ of $z_{0}$ where $u$ has a Holder continuous representative on $O_{z_{0}}\cap Q$.
 By  Lemma 3.3 of \cite{sersver02}, there exists a universal constant $c_{0}>0$ and a sequence of numbers  $R_{k}\in ]0,1[$ such that $R_{k}\rightarrow0$ as $k\rightarrow+\infty$ and
\begin{equation}\label{singcondition}
A(z_0,R_k;u)=\sup_{t_0-R_{k}^{2}\leqslant s\leqslant t_0}\frac{1}{R_{k}}\int\limits_{B(x_{0},R_k)}|u(x,s)|^2dx\geqslant c_0
\end{equation}
for any $k=1,2,\ldots$.\\
For each $\Omega=\omega\times ]a,b[$, where $\omega\Subset\mathbb{R}^{3}$ and $-\infty<a<b\leqslant 0$, we choose a large $k_{0}=k_{0}(\Omega)\geqslant 1$ so that for all $k\geqslant k_{0}$ we have for $(x,t)\in \Omega$:
$$x_{0}+R_k x\in B({2}/{3}),$$
and
$$t_0+R_{k}^{2}t\in ]-{({2}/{3}})^{2},0[.$$
Given such an $\Omega$ we perform the Navier-Stokes scaling as follows:
$$u^k(x,t):=R_ku(x_0+R_k x,t_0+R_{k}^{2}t)$$
and 
$$p^{k}(x,t):=R_{k}^2p(x_0+R_{k}x,t_{0}+R_{k}^{2}t).$$ 
As done in \cite{Phuc} we may decompose the pressure
\begin{equation}\label{presdecomp}
p=\tilde{p}+h,
\end{equation}
where $h$ is harmonic in $B$, and $\tilde{p}:=R_{i}R_{j}[(u_1u_j)\chi_{B}].$
We also have the decomposition for the rescaled pressure
\begin{equation}\label{rescaledpresdecomp}
p^k=\tilde{p^k}+h^k.
\end{equation}
Here,
$$\tilde{p^k}(x,t):=R_{k}^{2}\tilde{p}(x_0+R_{k}x,t_{0}+R_{k}^2 t),\,\,\,\,h^{k}(x,t):=R_{k}^{2}\tilde{h}(x_0+R_{k}x,t_{0}+R_{k}^2 t)$$
for any $(x,t)\in\Omega$ and $k\geqslant k_{0}(\Omega).$
Now under the hypothesis of Theorem \ref{Interiorreg},  Lemma \ref{usuitinter} gives that $(u,p)$ is a suitable weak solution to the Navier Stokes equations in $Q({5}/{6})$. It is not difficult to see that this implies $(u^k,p^k)$ is a suitable weak solution to the Navier Stokes equations in $\Omega$.\\

Now  define the relevant following scale invariant functional ($0<r<1$)
\begin{equation}\label{energy1}
A(z_{0},r;u):= \sup_{t_{0}-r^2\leqslant t\leqslant t_0}r^{-1}\int\limits_{B(x_{0},r)}|u(x,t)|^2dx,
\end{equation}
\begin{equation}\label{dissipation}
B(z_{0},r;u):=r^{-1}\int\limits_{Q(z_{0},r)}|\nabla u(x,t)|^2dxdt,
\end{equation}
\begin{equation}\label{uweak}
C_{\infty}(z_{0},r;u):= r^{-2}\int\limits_{t_{0}-r^2}^{t_0}\|u\|_{L^{3,\infty}(B(x_0,r))}^3 dt,
\end{equation}
\begin{equation}\label{weakpressure}
D_{\infty}(z_0,r;p):=r^{-2}\int\limits_{t_0-r^2}^{t_0}\|p\|_{L^{\frac{3}{2},\infty}(B(x_{0},r))}^{\frac{3}{2}}dt.
\end{equation}
 Now let us state a new estimate, which may be of independent interest, that we use heavily in Theorem \ref{Interiorreg} and Theorem \ref{Barkerboundaryreg}.
The proof of this is contained in the Appendix.
\begin{lemma}\label{energyboundLorentz}
Let $(u,p)$ be a suitable weak solution in $Q(z_{0},1)$. Then for $0<r<1$ the following  holds (c is some universal constant):
\begin{equation}\label{energyweaknorms1}
A(z_{0},{r}/{2};u)+ B(z_0,{r}/{2};u)\leqslant c( C_{\infty}(z_{0},r;u)^{\frac{4}{3}}+C_{\infty}(z_{0},r;u)^{\frac{2}{3}}+D_{\infty}(z_{0},r;p)^{\frac{4}{3}}).
\end{equation}
\end{lemma}
Now, let us make more explicit the role that this estimate plays by means of a Proposition. 
\begin{pro}\label{interioruniformestBarker}
The rescaled velocity and pressure the following uniform estimates  for $k\geqslant k_0(\Omega)$:
\begin{equation}\label{interiorBarkerpkuniform}
\|p^k\|_{L_{\frac{3}{2}}(a,b;L^{\frac{3}{2},\frac{q}{2}}(\omega)}\leqslant C(\Omega)[\|p\|_{L_{\frac{3}{2}}(Q)}+\|u\|_{L_{\infty}(-1,0;L^{3,q}(B)}^2],
\end{equation} 
\begin{equation}\label{rescaledharmonicest}
\int\limits_{a}^{b}\sup_{\omega}|h^k(x,t)|^{\frac 3 2}dt\leqslant
CR_k (\|p\|_{L_{\frac{3}{2}}(Q)}^{\frac{3}{2}}+
\|u\|_{L_{\infty}(-1,0;L^{3,q}(B)}^{3}),
\end{equation}
\begin{equation}\label{interiorBarkerukLorentz}
\|u^k(\cdot,t)\|_{L^{3,q}(\omega)}\leqslant\|u\|_{L_{\infty}(-1,0;L^{3,q}(B))}
\end{equation}
and
\begin{equation}\label{interiorBarkerukenergy}
\|u^k\|_{L_{\infty}(a,b;L_2(\omega)}+\|\nabla u^k\|_{L_{2}(a,b;L_{2}(\omega)}\leqslant C(\Omega,\|p\|_{L_{\frac{3}{2}}(Q)},\|u\|_{L_{\infty}(-1,0;L^{3,q}(B))})
\end{equation}
for all sufficiently large $k$ depending only on $\Omega$.
\end{pro}
\textbf{Proof}\\
It is clear  (\ref{interiorBarkerukLorentz}) follows from Lemma \ref{usuitinter} along with the easily seen property that if $u_{0}\in  L^{3,q}(\Omega)$ and $u_{0,\lambda}(x)= \lambda u_{0}(\lambda x)$ implies $\|u_{0,\lambda}\|_{L^{3,q}(\Omega/\lambda)}=\|u_{0}\|_{L^{3,q}(\Omega)}.$ Such spaces are called critical spaces for the Navier-Stokes equations.
Clearly, (\ref{interiorBarkerukenergy}) follows from Lemma \ref{energyboundLorentz} having established (\ref{interiorBarkerpkuniform})-(\ref{interiorBarkerukLorentz}).
Hence, we focus only on proving (\ref{interiorBarkerpkuniform})-(\ref{rescaledharmonicest}).
Note that we use the notation from (\ref{presdecomp})-(\ref{rescaledpresdecomp}).
  Lemma \ref{usuitinter} and inverse Navier-Stokes scaling gives us the following for $t\in ]a,b[$:
\begin{equation}\label{RescaledRieszest}
\|\tilde{p^{k}}(\cdot,t)\|_{L^{\frac{3}{2},\frac{q}{2}}(\omega)}\leqslant \|\tilde{p}(\cdot,t_0+R_{k}^{2}t)\|_{L^{\frac{3}{2},\frac{q}{2}}(B(\frac{2}{3}))}\leqslant C\|u\|_{L_{\infty}(-1,0;L^{3,q}(B))}^{2}.
\end{equation}
Here, we used the well known fact that Calderon Zygmund singular integral operators are bounded linear operators on $L^{p,q}$ with $p\in ]1,\infty[$ and $q\in ]0,\infty]$.
Using (\ref{RescaledRieszest}) we see it is sufficient to only prove (\ref{rescaledharmonicest}).
 Since $h$ is a harmonic function in $B$ for a.a $t\in ]-1,0[$, we have
\begin{equation}\label{harmonicpartest}
\|h\|_{L_{\frac{3}{2}}(-1,0;L_{\infty}(B(\frac{3}{4}))}\leqslant C
\|h\|_{L_{\frac{3}{2}}(-1,0; L_{1}(B(\frac{3}{4}))}\leqslant C(\|\tilde{p}\|_{L_{1,\frac{3}{2}}(Q)}+\|p\|_{L_{\frac 3 2}(Q)})\leqslant$$$$\leqslant C (\|\tilde{p}\|_{L_{\frac 3 2}(-1,0;L^{\frac 3 2,\frac q 2}(B))}+\|p\|_{L_{\frac 3 2}(Q)})\leq c(\|p\|_{L_{\frac 3 2}(Q)}+\|u\|_{L_{\infty}(-1,0; L^{3,q}(B))}^2).
\end{equation}
In the above, we used what is known as O'Neil inequality. Namely
if $0<p,q,r\leqslant\infty$, $0<s_1,s_2\leqslant\infty$,$\frac{1}{p}+\frac{1}{q}=\frac{1}{r}$ and $\frac{1}{s_1}+\frac{1}{s_2}=\frac{1}{s}$ then:
\begin{equation}\label{Holderverygeneral}
\|fg\|_{L^{r,s}(\Omega)}\leqslant C(p,q,s_1,s_2)\|f\|_{L^{p,s_1}(\Omega)}\|g\|_{L^{q,s_2}(\Omega)}.
\end{equation}
Using (\ref{harmonicpartest}), one infers 
\begin{equation}\label{rescaledharmonicest1}
\int\limits_{a}^{b}\sup_{\omega}|h^k(x,t)|^{\frac{3}{2}}dt\leqslant R_{k}\int\limits_{-({3}/{4})^{2}}^{0}\sup_{\omega}|h(x_0+R_k x,s)|^{\frac{3}{2}}ds\leqslant$$$$\leqslant R_{k}\|h\|_{L_{\frac{3}{2}}(-1,0;L_{\infty}(B(\frac{3}{4}))}^{\frac{3}{2}}\leqslant CR_k (\|p\|_{L_{\frac{3}{2}}(Q)}^{\frac{3}{2}}+
\|u\|_{L_{\infty}(-1,0;L^{3,q}(B)}^{3}).
\end{equation}
$\Box$\\
The remaining part of the proof of Theorem \ref{Interiorreg} is proven in more or less an identical way to the method in \cite{Phuc}. $\Box$
\setcounter{equation}{0}
\section{Local regularity near the flat part of the boundary in nonendpoint borderline Lorentz spaces}
Define the mixed Sobolev space:
$$ W^{2,1}_{m,n}(Q_{T})=\{ v\in L_{m,n}(Q_{T}): \|v\|_{L_{m,n}(Q_{T})}+$$$$+\|\nabla v\|_{L_{m,n}(Q_{T})}+\|\nabla^{2} v\|_{L_{m,n}(Q_{T})}+\|\partial_{t} v\|_{L_{m,n}(Q_{T})}<\infty\}.$$

First we show that the assumptions of Theorem \ref{Barkerboundaryreg} immediately give spatial smoothing and a local energy inequality.
\begin{pro}\label{localenergyboundarysmoothing}
Suppose $(v,p)$ satisfy the hypotheses of Theorem \ref{Barkerboundaryreg}.
Then for $\tau\in ]0,1[$ we have
\begin{equation}\label{uhigherreg}
v\in L_{4}(Q^{+}(2))\cap W^{2,1}_{\frac{9}{8},\frac{3}{2}}(Q^{+}(2\tau))\cap W^{2,1}_{\frac{4}{3}}(Q^{+}(2\tau)),
\end{equation}
\begin{equation}\label{phigherreg}
p\in W^{1,0}_{\frac{9}{8},\frac{3}{2}}(Q^{+}(2\tau))\cap W^{1,0}_{\frac{4}{3}}(Q^{+}(2\tau)). 
\end{equation}
In addition, the inequality
\begin{equation}\label{lorentzbddpointwisehalf}
\|v(\cdot,t)\|_{L^{3,q}(B^+(2\tau)}\leqslant \|v\|_{L_{\infty}(-(2^2),0; L^{3,q}(B^+(2))}
\end{equation}
holds for all $t\in ]-(2\tau)^2,0[$, and the function
$$ t\rightarrow \int\limits_{B(2\tau)}v(x,t)w(x) dx$$
is continuous on $[-(2\tau)^2,0]$ for any $w\in L^{\frac{3}{2},\frac{q}{q-1}}(B^+(2\tau)).$\\
Moreover, for a.a $t\in ]-1,0[$ and for all non negative functions $\phi\in C_{0}^{\infty}(\mathbb{R}^4)$, vanishing in a neighbourhood of the parabolic boundary $\partial^{'}Q(2)$ of $Q(2)$, $v$ and $p$ satisfy the inequality
\begin{equation}\label{localenergyinequalitytheo}
\int\limits_{B^+(2)}\phi(x,t)|v(x,t)|^2dx+2\int\limits_{B^{+}(2)\times ]-2,t[}\phi |\nabla v|^2 dxdt^{'}\leqslant$$$$\leqslant
\int\limits_{B^{+}(2)\times ]-2,t[}[|v|^2(\partial_{t}\phi+\Delta\phi)+v.\nabla\phi(|v|^2+2p)] dxdt^{'}.
\end{equation}
\end{pro}
\textbf{Proof}\\
First we remark that (\ref{localenergyinequalitytheo}) and (\ref{lorentzbddpointwisehalf}) is a consequence of (\ref{uhigherreg})-(\ref{phigherreg}), together we arguments used in the previous section on the interior case, so we focus on only proving (\ref{uhigherreg})-(\ref{phigherreg}). We proceed in a slightly different way to \cite{S2005}. Indeed we instead use the results regarding local boundary regularity for Stokes equation developed by Seregin in \cite{Ser09}. 
Under the assumptions of Theorem \ref{Barkerboundaryreg} we clearly have that
\begin{equation}\label{utheo1dif}
v\in W^{1,0}_{\frac{9}{8},\frac{3}{2}}(Q^{+}(2))\cap W^{1,0}_{\frac{4}{3}}(Q^{+}(2)),
\end{equation}
\begin{equation}\label{ptheo1dif}
p\in L_{\frac{9}{8},\frac{3}{2}}(Q^{+}(2))\cap L_{\frac{4}{3}}(Q^{+}(2)).
\end{equation}
Clearly from Sobolev embedding we have $v\in L_{6,2}(Q^{+}(2))$.
By a well known characterisation of Lorentz spaces:
\begin{equation}\label{Lorentzcharacter}
 L_{4}(B^+(2))=(L^{3,q}(B^+(2)),L_{6}(B^+(2)))_{\frac{1}{2},4}.
 \end{equation}
The proof of this and notation used in (\ref{Lorentzcharacter}) can be found in \cite{Bergh}, for example.
 By well known properties of interpolation spaces ( see, for example, section 3.5 of \cite{Bergh}):
$$\|v(\cdot,t)\|_{L_{4}(B^+(2))}\leqslant C\|v(\cdot,t)\|_{L_{3,q}(B^{+}(2))}^{\frac{1}{2}}\|v(\cdot,t)\|_{L_{6}(B^{+}(2))}^{\frac{1}{2}}.$$
Thus,
\begin{equation}\label{vL4}
\|v\|_{L_{4}(Q^{+}(2))}\leqslant C\|v\|_{L_{\infty}(-4,0;B^{+}(2))}^{\frac{1}{2}}\|v\|_{W^{1,0}_{2}(Q^+(2))}^{\frac{1}{2}}.
\end{equation}
By Holder's inequality we obtain
\begin{equation}\label{nonlinest4/3}
\|v.\nabla v\|_{L_{\frac{4}{3}}(Q^+(2))}\leqslant C\|v\|_{L_{\infty}(-4,0;B^{+}(2))}^{\frac{1}{2}}\|v\|_{W^{1,0}_{2}(Q^+(2))}^{\frac{3}{2}}.
\end{equation}
It is well known that the assumption $v\in W^{1,0}_{2}(Q^+(2))\cap L_{2,\infty}(Q^{+}(2))$ implies:
\begin{equation}\label{nonlinest9/8}
v.\nabla v \in L_{\frac{9}{8},\frac{3}{2}}(Q^+(2)). 
\end{equation} 
From here (\ref{utheo1dif})-(\ref{ptheo1dif}) combined with (\ref{nonlinest4/3})-(\ref{nonlinest9/8}) are enough to obtain the higher regularity.
This is seen from the local boundary regularity for the Stokes system namely Lemma 1.1 of \cite{Ser09}. $\Box$
\subsection{Pressure estimates}
It is first necessary to recap some appropriate definitions and a Lemma that will be necessary for our investigation.
The first definition was given by Seregin in \cite{S3}. Later on in \cite{Ser09}, Seregin gave a more general definition of suitable weak solution near flat part of boundary, but the one stated below is sufficient for our purposes.
\begin{definition}\label{suitableweakflatpart}

 A pair of functions $v$ and $p$ is called a suitable weak solution to the Navier Stokes equations in $Q^{+}(z_0,R)$ near $\Gamma(x_{0},R)\times [t_{0}-R^2,t_{0}]$ if they satisfy the following conditions.
They have the differentiability properties
\begin{equation}\label{vdifsuithalf}
v\in L_{\infty}(t_{0}-R^{2},t_{0}; L_{2}(B^{+}(x_{0},R))\cap W^{1,0}_{2}(Q^{+}(z_{0},R))\cap W^{2,1}_{\frac{9}{8},\frac{3}{2}}(Q^{+}(z_0,R)),
\end{equation}
\begin{equation}\label{pdifsuithalf}
p\in W^{1,0}_{\frac{9}{8},\frac{3}{2}}(Q^{+}(z_{0},R)).
\end{equation}
The pair $v$ and $p$ satisfies the Naiver Stokes equations a.e in $Q^{+}(z_{0},R)$ and the boundary condition
\begin{equation}\label{vsuitzerotrace}
v(x,t)=0,\,\,\,\,x_{3}=x_{03}\,\,and\,\, t_{0}-R^{2}< t< t_{0}.
\end{equation}
For a.a $t\in ]t_{0}-R^2,t_0[$ and for all non negative cut-off functions $\phi\in C^{\infty}_{0}(\mathbb{R}^4)$ vanishing in a neighbourhood of the parabolic boundary
$$\partial^{'}Q(z_0,R)= B(x_0,R)\times \{t=t_0-R^2\}\cup \partial B(x_0,R)\times [t_{0}-R^2,t_0]$$
of the cylinder $Q(z_{0},R)$, $v$ and $p$ satisfy the local energy inequality
\begin{equation}\label{localenergyinequalityflat}
\int\limits_{B^+(x_0,R)}\phi(x,t)|v(x,t)|^2dx+2\int\limits_{B^{+}(x_0,R)\times ]t_{0}-R^2,t[}\phi |\nabla v|^2 dxdt^{'}\leqslant$$$$\leqslant
\int\limits_{B^{+}(x_{0},R)\times ]t_{0}-R^2,t[}[|v|^2(\partial_{t}\phi+\Delta\phi)+v.\nabla\phi(|v|^2+2p)] dxdt^{'}.
\end{equation}
\end{definition} 
 Before stating and proving a certain Lemma let us introduce some notation.
 Various mean values of integrable functions are denoted as follows
$$[p]_{\Omega}=\frac{1}{|\Omega|}\int\limits_{\Omega}p(x,t)dx,$$
$$(v)_{\omega}=\frac{1}{\omega}\int\limits_{\omega}v dz.$$
Take $q\in[3,\infty]$ and introduce the following scale invariant quantities
\begin{equation}\label{vLorentzq}
C_{q}(z_{0},r;v):=\frac{1}{r^2}\int\limits_{t_0-r^2}^{t_0}\|v(\cdot,t)\|_{L^{3,q}(B(x_{0},r))}^3 dt
\end{equation}
\begin{equation}\label{pLorentzq}
D_{q}(z_{0},r;p):=\frac{1}{r^2}\int\limits_{t_0-r^2}^{t_0}\|p(\cdot, t)-[p]_{B(x_{0},r)}(t)\|_{L^{\frac{3}{2},\frac{q}{2}}(B(x_{0},r))}^{\frac{3}{2}} dt 
\end{equation}
The following generalises Lemma 2.1 of \cite{S2005}. Let us Remark that the Lemma and proof is very similar to Lemma 3.1 in \cite{ser01}. As it is slightly different we outline a proof for convenience of reader. Here it is:
\begin{lemma}\label{interiorpressureLorentzdecay}
Let $v\in L_3(Q(z_0,R))$ and $p\in L_{\frac{3}{2}}(Q(z_0,R))$ satisfy the Navier Stokes equations in the sense of distributions. Then, for $0<r\leqslant \rho\leqslant R$, we have
\begin{equation}\label{Lorentzpresinteriorest}
D_{q}(z_0,r;p)\leq c[\left(\frac{r}{\rho}\right)^{\frac{5}{2}}D_{q}(z_0,\rho;p)+\left(\frac{\rho}{r}\right)^{2}C_{q}(z_0,\rho;v)].
\end{equation}
\end{lemma}
This enables us to prove a certain generalisation to Proposition 2.5 in \cite{S2005}.
Namely the following.
\begin{pro}\label{presdecayinteriormainprop}
Assume all the criteria of Lemma \ref{interiorpressureLorentzdecay} are fulfilled. And let, in addition, for $q\in [3,\infty[$, that
\begin{equation}\label{interLorentzbddpro}
\|v\|_{L_{\infty}(t_0-R^2,t_0;L^{3,q}(B(x_0,R))}\leqslant L<\infty.
\end{equation}
Then, for any $\gamma\in ]0,1[$, there exists a constant $c_{1}(L,\gamma)$ such that, for $0<r\leqslant R$, we have
\begin{equation}\label{presdecayinterformula}
D_{q}(z_0,r;p)\leqslant c_1\left[\left(\frac{r}{R}\right)^{\frac{5}{2}\gamma}D_{q}(z_0,R;p)+1\right].
\end{equation}
\end{pro}
\textbf{Proof}\\
Clearly from (\ref{interiorpressureLorentzdecay}) we have
\begin{equation}\label{presinteriterate}
D_{q}(z_0,\tau^{k+1}R;p)\leqslant c\left[\tau^{\frac{5}{2}}D_{q}(z_0,\tau^{k}R;p)+\frac{L^{3}}{\tau^{2}}\right]
\end{equation}
for any $0<\tau<1$. We choose $\tau$ such that $c\tau^{\frac{5}{2}(1-\gamma)}\leq 1$. The proof then follows immediately from iterating (\ref{presinteriterate}).
$\Box.$\\
Before introducing the relevant statements we are required to introduce further notation. Firstly, quantities already defined but with  $'+'$ superscript denote integration over half-balls.
For example,
$$A^{+}(z_{0},r;u):= \sup_{t_{0}-r^2\leqslant t\leqslant t_0}r^{-1}\int\limits_{B^{+}(x_{0},r)}|u(x,t)|^2dx.$$

 Now define the following:
\begin{equation}\label{halfpressure}
D_{1}^{+}(z_{0},r;p):=\frac{1}{r^{\frac{3}{2}}}\int\limits_{t_0-r^2}^{t_0}(\int\limits_{B^{+}(x_{0},r)}|\nabla p|^{\frac{9}{8}}dx)^{\frac{4}{3}} dt.
\end{equation}
The following Lemma was proven in \cite{S3} and is Lemma 7.2 there. We state it without proof.
\begin{lemma}\label{seregin2002}
Consider a pair of functions $v$ and $p$  that is a suitable weak solution to the Navier Stokes equations in $Q^{+}(z_0,R)$ near $\Gamma(x_{0},R)\times [t_{0}-R^2,t_{0}].$ 
Then for any $0<r\leqslant\rho\leqslant R$, we have
\begin{equation}\label{sereginpressuredecayhalf}
D_{1}^+(z_0,r;p)\leqslant c\{\left(\frac{r}{\rho}\right)^{2}[ D_{1}^+(z_0,\rho;p)+(B^{+}(z_0,\rho;v))^{\frac{3}{4}}]+$$$$+\left(\frac{\rho}{r}\right)^{\frac{3}{2}}(A^+(z_0,\rho;v))^{\frac{1}{2}}B^+(z_0,\rho;v)\}.
\end{equation}
\end{lemma}
Next we state a new estimate. The statement and proof in the Appendix is for the interior case, but there is no distinction in the proof for the flat part of the boundary. It is very similar to Lemma \ref{energyboundLorentz}, and is a key part that allows us to work in the context of nonendpoint critical Lorentz spaces.

\begin{lemma}\label{Lorentzenergygradp2}
Let $(u,p)$ be a suitable weak solution in $Q^{+}(z_{0},1)$ near $\Gamma(x_0,1)\times [t_0-1,t_0].$ 
Then for $0<r<1$ the following  holds (c is some universal constant):
\begin{equation}\label{energyweaknorms1half}
A^{+}(z_{0},{r}/{2};u)+ B^{+}(z_0,{r}/{2};u)\leqslant c( C^{+}_{\infty}(z_{0},r;u)^{\frac{4}{3}}+C^{+}_{\infty}(z_{0},r;u)^{\frac{2}{3}}+$$$$+D^{+}_{1}(z_{0},r;p)^{\frac{2}{3}}C_{\infty}^{+}(z_{0},r;u)^{\frac{1}{3}}).
\end{equation}
\end{lemma}
Now we can prove a generalisation of Proposition 2.6 in \cite{S2005}.
\begin{pro}\label{presdecayhalf}
Let $q\in [3,\infty[.$
Assume that $(v,p)$ is a suitable weak solution in $Q^{+}(z_0,R)$ near $\Gamma(x_{0},R)\times [t_{0}-R^2,t_{0}]$. Suppose, in addition, 
\begin{equation}\label{Lorentzbddhalf}
\|v\|_{L_{\infty}(t_0-R^2,t_0; L^{3,q}(B(z_0,R))}\leqslant L<\infty.
\end{equation}
Then, for any $\gamma\in ]0,1[$, there exists a constant $c_{2}$ depending on $\gamma$ and $L$ only such that, for $0<r\leqslant R$, we have
\begin{equation}\label{presdecayhalfformula}
D_{1}^{+}(z_0,r;p)\leqslant c_{2}\left[\left(\frac{r}{R}\right)^{2\gamma}D_{1}^{+}(z_0,R;p)+1\right].
\end{equation}
\end{pro}
\textbf{Proof}\\
Let $\rho\leqslant\frac{R}{2}$. Then by Lemma \ref{Lorentzenergygradp2} and (\ref{Lorentzbddhalf}) we have
$$ A^{+}(z_0,\rho;v)+B^{+}(z_0,\rho;v)\leq c[L^{4}+L^2+ (D_{1}^{+}(z_{0},2\rho;p))^{\frac{2}{3}}L].$$
By O'Neils inequality for $0\leqslant s\leqslant R$:
$$\int\limits_{B(x_{0},s)} |v(x,t)|^2 dx\leqslant C|B(0,s)|^{\frac{1}{3}}\||v|^2(\cdot,t)\|_{L^{\frac{3}{2},\infty}(B(x_{0},s))}.$$
Thus, it is clear that
\begin{equation}\label{AboundL}
A^+(z_0,s;v)\leqslant C L^2.
\end{equation}
Thus using these facts with by Lemma \ref{seregin2002}, we see that for any $0<r\leqslant\rho\leqslant\frac{R}{2}$:
$$ D_{1}^{+}(z_0,r;p)\leq c\left(\frac{r}{\rho}\right)^{2}[ D_{1}^+(z_0,\rho;p)+((L^{4}+L^2+ (D_{1}^{+}(z_{0},2\rho;p))^{\frac{2}{3}}L)^{\frac{3}{4}}]+$$$$+\left(\frac{\rho}{r}\right)^{\frac{3}{2}}L [L^{4}+L^2+ (D_{1}^{+}(z_{0},2\rho;p))^{\frac{2}{3}}L]\}.
$$
Then by Young's inequality, it is not so difficult to see for $0<r\leqslant \rho\leqslant \frac{R}{2}$:
$$ D_{1}^+(z_0,r;p)\leqslant c(L)\left[\left(\frac{r}{\rho}\right)^{2}(D_{1}(z_0,2\rho;p)+1)+\left(\frac{\rho}{r}\right)^{\frac{17}{2}}\right].$$
Clearly this implies that for $0<r\leqslant\rho\leqslant R$ that
$$ D_{1}^+(z_0,r;p)\leqslant c(L)\left[\left(\frac{r}{\rho}\right)^{2}(D_{1}(z_0,\rho;p)+1)+\left(\frac{\rho}{r}\right)^{\frac{17}{2}}\right].$$
The conclusion now follows by identical reasoning as Proposition \ref{presdecayinteriormainprop}. $\Box$
\subsection{Proof of Theorem 1.3}
We let $$L=\|v\|_{L_{\infty}(-4,0; L^{3,q}(B^+(2)))}<\infty.$$
Let us describe the rescaling procedure taken from \cite{S2005}. Assume Theorem \ref{Barkerboundaryreg} is false. Let $z_{0}\in \bar{Q}^{+}({1}/{2})$ be a singular point.
 Then we know that all conclusions of Proposition \ref{localenergyboundarysmoothing} hold. Hence, $(v,p)$ form  a suitable weak solution to the Navier Stokes equations in $Q^{+}(1)$ near $\Gamma(0,1)\times[-1,0].$ Furthermore for every $t\in [-1,0]$
 \begin{equation}\label{pointwiseLorentzbdd}
 \|v(\cdot,t)\|_{L^{3,q}(B^+(1))}\leqslant L.
 \end{equation}
 From Theorem \ref{Interiorreg} we can only have $z_{0}$ lies on the boundary $\bar{\Gamma}(0,{1}/{2}).$
 Without loss of generality (translation invariance), we assume that $z_{0}=0$ and that all conclusions of Proposition \ref{localenergyboundarysmoothing} (together with (\ref{pointwiseLorentzbdd})) hold on the slightly smaller domain
 $Q^{+}({1}/{2})$.\\
 As a consequence of Lemma 3.3, proven in \cite{S2005}, there exists a decreasing sequence $R_{k}<\frac{1}{2}$ tending to zero, together with a universal constant $\epsilon_{3}$, such that for $k=1,2\ldots$
 \begin{equation}\label{epsilonreghalf}
 \frac{1}{R_{k}^{2}}\int\limits_{Q^{+}(R_k)}|v|^3 dz\geqslant \epsilon_3.
  \end{equation}
  Extending functions $(v,p)$ outside $Q^{+}({1}/{2})$ to zero, for $ (y,s)\in\mathbb{R}^{3}_{+}\times ]-\infty,0[$ define the rescaled functions
  $$ u^k(y,s)= R_{k}v(R_ky,R_k^2s),\,\,\,\, p^k(y,s)=R_{k}^{2}p(R_{k}y,R_{k}^{2}s).$$
  Next we claim the following properties in the limit:
\begin{pro}\label{uklimhalf}
There exists a subsequence of $(u^k,p^k)$, still denoted by $(u^k,p^k)$, and a pair of functions $(u_{\infty},p_{\infty})$ with $\rm{div}\, u_{\infty}=0$ in $\mathbb{R}^{3}_{+}\times ]-\infty,0[$, such that 
\begin{equation}\label{lorentzhalfweakstar}
u^{k}\stackrel{*}{\rightharpoonup} u_{\infty}\,\,\, \rm{in}\,\,L_{\infty}(-\infty,0; L^{3,q}(\mathbb{R}^{3}_{+}).
\end{equation}
Moreover for any $a>0$
\begin{equation}\label{locallyufinitehalf}
|u_{\infty}|^{2},\,\nabla u_{\infty}\in L_{2}(Q^{+}(a)).
\end{equation}
Additionally,
\begin{equation}\label{locallypfinitehalf}
(u_{\infty},p_{\infty})\in W^{2,1}_{\frac{4}{3}}(Q^+(a))\cap W^{2,1}_{\frac{9}{8},\frac{3}{2}}(Q^+(a))\times W^{1,0}_{\frac{4}{3}}(Q^+(a))\cap W^{1,0}_{\frac{9}{8},\frac{3}{2}}(Q^+(a)),
\end{equation}
and $(u_{\infty},p_{\infty})$ forms a suitable weak solution to the Navier-Stokes equations in $Q^+(a)$ near $\Gamma(a)\times [-a^2,0]$.
One also has that
\begin{equation}\label{convergcontin}
u^{k}\rightarrow u_{\infty}\,\,\rm{in}\,\, C([-a^{2},0]; L_{s}(B^{+}(a)),\,\ u_{\infty}(x,0)=0
\end{equation}
for any $s\in ]1,3[$ and for a.a $x\in\mathbb{R}^{3}_{+}$.
Furthermore, $u_{\infty}$ satisfies the lower bound
\begin{equation}\label{uinfinitylowerbound}
\int\limits_{Q^+}|u_{\infty}|^3 dz\geqslant \epsilon_{3}.
\end{equation}
\end{pro}
\textbf{Proof}\\
Fix $a>0$ and let $k(a)$ be such that
\begin{equation}\label{Rkcondition}
R_{k}a<\frac{1}{8}
\end{equation}
for all $k\geqslant k(a)$. It is clear that $(\ref{lorentzhalfweakstar})$ follows from identical reasons as discussed in the interior case, and in addition
\begin{equation}\label{Lorantzbddlimit}
\|u_{\infty}\|_{L_{\infty}(-\infty,0; L^{3,q}(\mathbb{R}^{3}_{+}))}\leqslant \sup_{k}\|u_{k}\|_{L_{\infty}(-\infty,0; L^{3,q}(\mathbb{R}^{3}_{+}))}=L<\infty.
\end{equation}
 By inverse scaling and  the same reasons discussed in Proposition \ref{presdecayhalf} we have 
 \begin{equation}\label{presenergysequence}
 B^{+}(0,2a; u^k)+D_{1}^{+}(0,2a,p^k)=B^{+}(0,2aR_k; v)+D_{1}^{+}(0,2aR_k,p)\leqslant$$$$\leqslant c_0(L)[1+D_{1}^{+}(0,4aR_k,p)].
 \end{equation}
 By Proposition \ref{presdecayhalf}
 $$D_{1}^{+}(0,4aR_k,p)\leqslant c_{1}(L)[8aR_k D_{1}^{+}(0,{1}/{2},p)+1]\leqslant c_{1}(L)[ D_{1}^{+}(0,{1}/{2},p)+1].$$
 Thus,
 \begin{equation}\label{uniformbddboundary}
  B^{+}(0,2a; u^k)+D_{1}^{+}(0,2a,q^k)\leqslant c_2(L)[ D_{1}^{+}(0,{1}/{2},p)+1].
 \end{equation}
 Since, $$L_{4}(B^+(a))= (L_{3,q}(B^+(a)), L_{6}(B^+(a)))_{\frac{1}{2},4},$$
 we have the interpolative inequality
 $$\|u^k(\cdot,t)\|_{L_{4}(B^+(2a))}\leq C\|u^k(\cdot,t)\|_{L^{3,q}(B^+(2a))}^{\frac{1}{2}}\|u^k(\cdot,t)\|_{L_{6}(B^+(2a))}^{\frac{1}{2}}\leqslant$$$$\leqslant CL^{\frac{1}{2}}\|u^k(\cdot,t)\|_{L_{6}(B^+(2a))}^{\frac{1}{2}}.$$
 From here, it follows that 
 \begin{equation}\label{L4seqbdd}
 \|u^{k}\|_{L_{4}(Q^+(2a))}\leqslant c_{0}(L,a)( D_{1}^{+}(0,{1}/{2},p)+1)^{\frac{1}{4}}.
 \end{equation}
 Observe that
 \begin{equation}\label{nonlin4/3half}
 \|u^k.\nabla u^k\|_{L_{\frac{4}{3}}(Q^+(2a))}\leqslant  \|u^{k}\|_{L_{4}(Q^+(2a))} \|\nabla u^{k}\|_{L_{2}(Q^+(2a))}\leqslant$$$$\leqslant c_{1}(L,a)( D_{1}^{+}(0,{1}/{2},p)+1)^{\frac{3}{4}}.
 \end{equation}
 Using multiplicative inequalities it is not so hard to see
 \begin{equation}\label{nonlin9/83/2half}
  \|u^k.\nabla u^k\|_{L_{\frac 9 8,\frac{3}{2}}(Q^+(2a))}\leqslant C_{2}(a,L) (B^{+}(0, 2a,u^{k})^{\frac{2}{3}}\leqslant$$$$\leqslant C_{3}(a,L)( D_{1}^{+}(0,{1}/{2},p)+1)^{\frac{2}{3}}.
 \end{equation}
 One can then use identical arguments to those used in \cite{S2005} to show that
 \begin{equation}\label{upbddhalf4/3}
 \|u^{k}\|_{W^{2,1}_{\frac{4}{3}}(Q^{+}(a))}+\|\nabla p^{k}\|_{L_{\frac{4}{3}}(Q^{+}(a))}\leqslant c_{4}(a,L,D_{1}^{+}(0,{1}/{2},p)), 
 \end{equation}
 \begin{equation}\label{upbddhalf9/83/2}
 \|u^{k}\|_{W^{2,1}_{\frac{9}{8},\frac{3}{2}}(Q^{+}(a))}+\|\nabla p^{k}\|_{L_{\frac 9 8,\frac{3}{2}}(Q^{+}(a))}\leqslant c_{5}(a,L,D_{1}^{+}(0,{1}/{2},p)).
 \end{equation}
 Having obtained this, (\ref{lorentzhalfweakstar})-(\ref{locallypfinitehalf}) and (\ref{uinfinitylowerbound}) follow in the same way as presented in \cite{S2005}. Let us focus on (\ref{convergcontin}). 
 By (\ref{upbddhalf4/3}) and compactness of the embedding  $$W^{2,1}_{\frac{4}{3}}(Q^+(a))\hookrightarrow C([-a^2,0]; L_{\frac{4}{3}}(B^{+}(a)),$$
 we have, taking further subsequences and using further cantor diagonalisation if necessary, that for all $a>0$
 \begin{equation}\label{4/3continconverg}
 u^{k}\rightarrow u_{\infty}\,\, \rm{in} \,\,\ C([-a^2,0]; L_{\frac{4}{3}}(B^{+}(a)).
 \end{equation}
 For $s\in ]\frac{4}{3},3[$ there exists $\theta(s)\in ]0,1[$ such that
 $$L_{s}(B^{+}(a))=(L_{\frac{4}{3}}(B^{+}(a)), L^{3,q}(B^{+}(a))_{\theta(s),s}$$
 using this and fact that $u^k$ and $u$ are uniformly bounded in $L_{\infty}(-a^2,0; L^{3,q}(B^{+}(a))$ it is simple to infer, by interpolative inequalities, that
 $$u^{k}\rightarrow u \,\, \rm{in}\,\, C([-a^2,0]; L_{s}(B^{+}(a)).$$
 For any $y\in B^ {+}(a)$, using O'Neil's inequality and inverse scailing gives us
 $$\int\limits_{B^{+}(y,1)} |u_{\infty}(x,0)| dx\leqslant\int\limits_{B^{+}(y,1)}|u_{\infty}(x,0)-u^{k}(x,0)| dx+\int\limits_{B^{+}(y,1)}|u^{k}(x,0)| dx\leqslant$$$$\leqslant \|u_{\infty}-u^{k}\|_{C([-(a+1)^2,0]; L_{1}(B^+(a+1))}+ |B_{1}(0)|^{\frac{2}{3}}\|v(\cdot,0)\|_{L^{3,q}(B^{+}(R_{k}y, R_{k}))}.$$
 Then using Proposition \ref{localenergyboundarysmoothing} and the dominated convergence theorem, applied to the distribution function of $v(\cdot,0)$, we see that
 $$\int\limits_{B^{+}(1)}|u_{\infty}(x,0)| dx=0.$$
  $\Box$\\
 It is necessary to introduce the following notation. Let $i_3=(0,0,1)$.
 The following is a generalisation of a result of \cite{S2005} (Lemma 4.1 there) .
 Here it is.
 \begin{lemma}\label{presbddweakhalf}
 There exists a positive constant $c_{5}(L, D_{1}^{+}(0,{1}/{2},p))$ with the following property. Fix $h>0$ and $T>0$ arbitrarily, then
 \begin{equation}\label{presbddweakhalfformula}
 D_{q}(e_0,2h;p_{\infty})\leqslant c_{5}
 \end{equation}
 for any $e_{0}= (y_{0},s_{0})\in (\mathbb{R}^{3}_{+}+3hi_3)\times ]-T,0[.$
 \end{lemma}
 \textbf{Proof}\\
 Let $a$ be sufficiently large such that $Q(e_{0},2h)\in Q^{+}(a).$
 From Proposition \ref{uklimhalf} and Poincare inequality we have that there exists $C_{a}(t)\in L_{\frac{3}{2}}(-a^2,0)$ such that
 $$ p^k-[p^{k}]_{B^{+}(a)}(t)  \rightharpoonup p_{\infty}+C_{a}(t)\,\,\rm{in}\,\, L_{\frac{3}{2}}(-a^2,0;L^{\frac{3}{2},\frac{q}{2}}(B^+(a)).  $$
 It easily follows from this that
 $$p^k-[p^{k}]_{B(y_0,2h)}(t)  \rightharpoonup p_{\infty}-[p_{\infty}]_{B(y_0,2h)}(t)\,\,\rm{in}\,\, L_{\frac{3}{2}}(-a^2,0;L^{\frac{3}{2},\frac{q}{2}}(B^+(a)).  $$
 Thus it is clear that
 $$\lim\sup_{k\rightarrow\infty} D_{q}(e_{0},2h; p^{k})\geqslant D_{q}(e_{0},2h; p_{\infty}).$$
 So, it is sufficient to prove the following bound
 \begin{equation}\label{presbddlargek}
 D_{q}(e_0,2h;p^{k})\leqslant c_{5}(L, D_{1}^{+}(0,{1}/{2},p))
 \end{equation}
for all $k$ sufficiently large such that
\begin{equation}\label{coordslargek}
x_{0}^k= y_{0}R_{k}\in B^{+}({1}/{8}),\,\,\,\,t_{0}^{k}= s_{0}R_{k}^{2}>-({1}/{8})^{2}.
\end{equation}
By inverse scaling, we have
\begin{equation}\label{pressureinverse}
D_{q}(e_0,2h;p^k)= D_{q}(z_{0}^k, 2hR_{k};p),\,\,\,\, z_{0}^{k}=(x_{0}^{k},t_{0}^{k}).
\end{equation}
Notice that $d_{k}= x_{03}^k= \rm{dist}(x_{0}^k,\Gamma)\geqslant 3hR_{k}$. Furthermore, $Q(z_{0}^k,2hR_{k})\subset Q(z_{0}^k,d_{k})\subset Q^{+}({1}/{4}).$ Thus from Proposition \ref{presdecayinteriormainprop} we see that
\begin{equation}\label{presbddlemma1}
D_{q}(z_{0}^k, 2hR_{k};p)\leqslant c(L)\left[\left(\frac{2hR_k}{d_k}\right)^{\frac{5}{4}}D_{q}(z_{0}^k, d_k;p)+1\right]\leqslant$$$$\leqslant
c(L)\left[D_{q}(z_{0}^k, d_k;p)+1\right].
\end{equation}
Let $\bar{x}_{0}^{k}= (x_{01}^{k}, x_{02}^{k},0)$ and also $\bar{z}_{0}^{k}=(\bar{x}_{0}^{k}, t_{0}^{k})$.
It is clear that $Q(z_{0}^{k},d_k)\subset Q^{+}(\bar{z}_{0}^{k}, 2d_{k})$ and 
$$ D_{q}(z_{0}^k, d_k;p)\leqslant cD_{q}^{+}(\bar{z}_{0}^k, 2d_k;p).$$
Note that we have the following Poincare inequality:
$$D_{q}^{+}(\bar{z}_{0}^k, 2d_k;p)\leqslant D_{1}^{+}(\bar{z}_{0}^k, 2d_k;p).$$
Using this and (\ref{presbddlemma1}) we infer
\begin{equation}\label{presbddbygrad}
D_{q}(z_{0}^k, 2hR_{k};p)\leqslant c(L)[D_{1}^{+}(\bar{z}_{0}^k, 2d_k;p)+1].
\end{equation}
Clearly $2d_{k}<\frac{1}{4}$ so we have
$$Q^{+}(\bar{z}_{0}^{k}, 2d_{k})\subset Q^{+}(\bar{z}_{0}^{k},{1}/{4})\subset Q^{+}({1}/{2}).$$
Then one can apply Proposition \ref{presdecayhalf} to infer that
\begin{equation}\label{gradpresbddprop}
D_{1}^{+}(\bar{z}_{0}^k, 2d_k;p)\leqslant c(L)[8d_{k}D_{1}^{+}(\bar{z}_{0}^k, 2d_k;p)+1]\leqslant$$$$\leqslant c(L)[D_{1}^{+}(0, {1}/{2};p)+1] .
\end{equation}
Thus putting everything together gives
$$D_{q}(e_{0}, 2h;q^k)\leqslant c(L)[D_{1}^{+}(0, {1}/{2};p)+1].$$ This easily gives us the conclusion by initial remarks. $\Box$\\
We now state a new regularity criteria, which is convenient for proving Theorem \ref{Barkerboundaryreg}. It is a generalisation of the basic $\epsilon$-regularity criteria found in \cite{Lin} and \cite{Ladyser02}.
The proof is contained in the Appendix. Here is the statement.
\begin{theorem}\label{strengthendepsilonregularity}
Let $(u,p)$ be a suitable weak solution in $Q_{1}(0)$. Then there exists a universal constants $\epsilon_{0}$ and $c_{0k}$ (with $k=1,2\ldots$) with the following property.  Assume
\begin{equation}\label{smallnesscondition}
C_{\infty}(0,1;u)+D_{\infty}(0,1;p)<\epsilon_{0}.
\end{equation}
then for any natural number $k$, $\nabla ^{k-1}u$ is Holder continuous in $\Bar{Q}({1}/{2})$ and the following bound is valid:
\begin{equation}\label{spatialsmoothness}
\max_{\Bar{Q}({1}/{2})}|\nabla^{k-1}u(z)|<c_{ok}.
\end{equation}

\end{theorem}
Now we proceed to the proof of Theorem \ref{Barkerboundaryreg}. Fix $h\in ]0,1[$ arbitrarily and consider an arbitrary point
$$z_0\in (\mathbb{R}^{3}_{+}+3hi_3)\times ]-100,0[.$$  Following \cite{S2005}, we perform the same pressure decomposition as in the interior case as follows.
In the ball $B(x_0,2h)$ decompose the pressure
$$p_{\infty}=p_{\infty}^1+p_{\infty}^2$$
such that
\begin{equation}\label{pinfty1def}
p_{\infty}^1:= R_{i}R_{j}((u_{\infty})_i(u_{\infty})_j\chi_{B(x_0,2h)}).
\end{equation}
It is clear that
$$ \Delta p_{\infty}^2(\cdot,t)=0\,\,\,\,in\,\,\,\, B(x_{0},2h).$$
Thus for the same reasons as previously stated we have
\begin{equation}\label{pres1inftyrieszest}
\|p_{\infty}^1(\cdot,t)\|_{L^{\frac{3}{2},\frac{q}{2}}(B(x_{0},h)}\leqslant c\|u_{\infty}(\cdot,t)\|_{L^{3,q}(B(x_{0},h)}^{2}.
\end{equation}
It is not so difficult to show that for the harmonic part $p_{\infty}^{2}(\cdot,t)$ we have the following:
$$\sup_{x\in B(x_0,h)}|\nabla p_{\infty}^{2}(x,t)|\leqslant \frac{c}{h}(\frac{1}{h^3}\int\limits_{B(x_{0},2h)}|p_{\infty}^{2}(\cdot,t)-[p_{\infty}^{2}]_{B(x_0,2h)}(t)|dx).$$
Clearly, by O'Neil's inequality:
\begin{equation}\label{pres2inftyharmonicest}
\sup_{x\in B(x_0,h)}|\nabla p_{\infty}^{2}(x,t)|\leqslant \frac{c}{h}\left(\frac{1}{h^2}\|p_{\infty}^{2}(\cdot,t)-[p_{\infty}^{2}]_{B(x_0,2h)}(t)\|_{L^{\frac{3}{2},\frac{q}{2}}(B(x_0,2h))}\right).
\end{equation}
For any $0<\rho<1$ we can use (\ref{pres1inftyrieszest})-(\ref{pres2inftyharmonicest}) along with the Poincare inequality to infer
$$ D_{q}(z_0,h\rho; p_{\infty})\leqslant c[D_{q}(z_0,h\rho; p_{\infty}^{1})+D_{q}(z_0,h\rho; p_{\infty}^{2})]\leqslant$$$$\leqslant \frac{c}{(h\rho)^2}\int\limits_{t_0-4h^2}^{t_0}\|p_{\infty}^1(\cdot,t)\|_{L^{\frac{3}{2},\frac{q}{2}}(B(x_{0},2h))}^{\frac{3}{2}} dt+$$$$+ c(h\rho)^{\frac{5}{2}}\int\limits_{t_0-(h\rho)^2}^{t_0}\sup_{x\in B(x_0,h)}|\nabla p_{\infty}^{2}(x,t)|^{\frac{3}{2}}dt\leqslant$$$$\leqslant
c\left[\frac{1}{(h\rho)^2}\int\limits_{t_0-4h^2}^{t_0}\|u_{\infty}(\cdot,t)\|_{L^{3,q}(B(x_{0},h))}^{3} dt+\rho^{\frac{5}{2}}D_{q}(z_0,2h,q)\right].$$
By Lemma \ref{presbddweakhalf} ($c_{5}=c_{5}(L,D_{1}^{+}(0,{1}/{2},p))$) and the well known embeddings for Lorentz spaces we obtain:
$$C_{\infty}(z_0,h\rho;u_{\infty})+D_{\infty}(z_0,h\rho;p_{\infty})
\leqslant$$$$\leqslant c\left[\frac{1}{(h\rho)^2}\int\limits_{t_0-4h^2}^{t_0}\|u_{\infty}(\cdot,t)\|_{L^{3,q}(B(x_{0},h))}^{3} dt+\rho^{\frac{5}{2}}c_5\right]$$
for any $z_0\in(\mathbb{R}^{3}_{+}+3hi_3)\times ]-100,0[.$
Next, fix $\rho(L, D_{1}^+(0,\frac{1}{2};p)\in ]0,1[$ such that
$$ c\rho^{\frac{5}{2}}c_{5}<\frac{\epsilon_{0}}{2}.$$
Since $u_{\infty}\in L_{\infty}(-\infty,0; L^{3,q}(\mathbb{R}^{3}_{+}))$ and $q\in ]3,\infty[$ we have that for a.a $t$
$$\lim_{R\rightarrow\infty}\|u_{\infty}(\cdot,t)\|_{L^{3,q}(\mathbb{R}^{3}_{+}\setminus B(R))}=0.$$
Thus, by the dominated convergence theorem, we may find $R_{1}>100$ such that
$$\frac{c}{(h\rho)^2}\int\limits_{-200}^{0}\|u_{\infty}(\cdot,t)\|_{L^{3,q}(\mathbb{R}^{3}_{+}\setminus B(\frac{R_1}{4}))}^{3} dt<\frac{\epsilon_{0}}{2}.$$
Using this, it can be inferred that for any $z_{0}\in (\mathbb{R}^{3}_{+}+3hi_3)\times ]-100,0[$ such that $|x_0|>\frac{R_1}{2}$ we have
$$C_{\infty}(z_0,h\rho;u_{\infty})+D_{\infty}(z_0,h\rho;p_{\infty})
<\epsilon_{0}.$$
Now Theorem \ref{strengthendepsilonregularity} is applicable. Moreover we obtain that  any $z_{0}\in (\mathbb{R}^{3}_{+}+6hi_3)\setminus B(R_1)\times ]-50,0[$  and $k=1,2\ldots$ 
$$|\nabla^{k}u_{\infty}(z_0)|\leqslant\frac{c_{0k}}{(h\rho)^{k+1}}.$$
By the interior result of Theorem \ref{Interiorreg}, we obtain boundedness of  $\nabla^{k}u_{\infty}$ on the set $(\mathbb{R}^{3}_{+}+6hi_3)\cap B(R_1)\times ]-50,0[$.
 Define the vorticity $\omega_{\infty} = \nabla \wedge u_{\infty}$. Then on the set 
$(\mathbb{R}^{3}_{+}+6hi_3)\times ]-50,0[$, we have that there exists $M>0$ such that
$$|\omega_{\infty}|\leqslant M,$$
$$|\partial_{t}\omega_{\infty}-\Delta\omega_{\infty}|\leqslant M(|\omega_{\infty}|+|\nabla \omega_{\infty}|).$$
Furthermore,
$$\omega_{\infty}(\cdot,0)=0\,\,\,\,\rm{in}\,\,\,\,\mathbb{R}^{3}_{+}.$$
By applying the backward uniqueness theorem (Theorem 5 in \cite{ESS2003}) one deduces that
$$\omega_{\infty}=0\,\,\,\,\rm{in}\,\,\,\,(\mathbb{R}^3_{+}+6hi_{3})\times [-50,0].$$
Since $h$ is arbitrary we infer that 
$$\omega_{\infty}=0\,\,\,\,\rm{in}\,\,\,\,\mathbb{R}^3_{+}\times [-50,0].$$
Hence for a.a $t\in [-50,0]$, $u_{\infty}$ is a harmonic function, which satisfies the boundary condition $u_{\infty}(x,t)=0$ if $x_{3}=0$. But for a.a $t\in [-50,0]$, $L^{3,q}$ of $u_{\infty}$ over $\mathbb{R}^{3}_{+}$ is finite. By a Liouville Theorem we get for the same $t$ that $u_{\infty}(\cdot,t)=0$ in $\mathbb{R}^{3}_{+}$. This contradicts (\ref{uinfinitylowerbound}). $\Box$

\setcounter{equation}{0}
\section{Appendix: Epsilon regularity in weak Lebesgue spaces}

First, we begin with stating a well known algebraic Lemma, whose proof is omitted but found in \cite{Giusti}.
\begin{lemma}\label{iterativeGiusti}
Let $I(s)$ be a bounded non negative function in the interval $[R_{1},R_{2}]$. Assume that for every $s,\rho\in [R_{1},R_{2}]$ and $s<\rho$ we have
$$ I(s)\leqslant [A(\rho-s)^{-\alpha}+B(\rho-s)^{-\beta}+C]+\theta I(\rho)$$
with $A,B,C\geqslant 0$, $\alpha>\beta>0$ and $\theta\in [0,1[$.
Then there holds
$$I(R_{1})\leqslant c(\alpha,\theta)[A(R_{2}-R_{1})^{-\alpha}+B(R_{2}-R_{1})^{-\beta}+C].$$
\end{lemma}

\textbf{Proof of Lemma \ref{energyboundLorentz}}
Without loss of generally, consider $z_{0}$ to be the origin.
Let $0<\frac{r}{2}\leqslant s<\rho\leqslant r<1$. Let $\eta_{1}\in C_{0}^{\infty}(B({\rho}))$ such that $0\leqslant \eta_{1}\leqslant 1$ in $\mathbb{R}^{3}$ and $\eta_{1}=1$ on $B({s})$.
Furthermore for $|\alpha|\leq 2$:
$$|\nabla^{\alpha}\eta_{1}|\leqslant \frac{C}{(\rho-s)^\alpha}.$$
Let $\eta_{2}\in C_{0}^{\infty}(-\rho^2,\rho^2)$ such that $0\leqslant \eta_{2}\leqslant 1$ in $\mathbb{R}$ and $\eta_{1}=1$ on $[-s^2,s^2]$.
Furthermore :
$$|\eta_{1}^{'}|\leqslant \frac{C}{(\rho^{2}-s^2)}\leqslant \frac{C}{r(\rho-s)}\leqslant \frac{C}{(\rho-s)^2}.$$
Let $\phi(x,t):=\eta_{1}(t)\eta_{2}(x)$.
Hence:
\begin{equation}\label{gradienttest}
|\nabla\phi|\leqslant\frac{C}{\rho-s},
\end{equation}

\begin{equation}\label{seconddertest}
|\nabla^{2}\phi|\leqslant\frac{C}{(\rho-s)^{2}},
\end{equation}

\begin{equation}\label{timedertest}
|\phi_{t}|\leqslant\frac{C}{(\rho-s)^{2}}.
\end{equation}
From the local energy inequality we have that for a.a $t\in ]-1,0[$:
\begin{equation}\label{localenergyinequality}
\int\limits_{B}\phi(x,t)|u(x,t)|^2dx+2\int\limits_{B\times ]-1,t[}\phi|\nabla u|^{2}dxds$$$$\leqslant\int\limits_{B\times ]-1,t[}(|u|^2(\Delta\phi+\partial_{t}\phi)+u.\nabla\phi(|u|^2+2p )dxds.
\end{equation}
Let $I_{1}(s):=sA(0,s;u)$, $I_{2}(s):=sB(0,s;u)$ and $I(s)=I_1(s)+I_{2}(s)$.
\begin{equation}\label{energyest1}
I(s)\leqslant (1)+(2)+(3).
\end{equation}
Where,
\begin{equation}\label{(1)}
(1):=\int\limits_{B\times ]-1,t[}(|u|^2(\Delta\phi+\partial_{t}\phi)dxds,
\end{equation}
\begin{equation}\label{(2)}
(2):= \int\limits_{Q(\rho)}u.\nabla\phi|u|^2dxds,
\end{equation}
\begin{equation}\label{(3)}
(3):=\int\limits_{Q(\rho)}2u.\nabla\phi pdxds.
\end{equation}
Let us treat simplest integral $(1)$ first.
By O'Neil's inequality in space and then Holder in time:
$$(1)\leqslant\int\limits_{-\rho^{2}}^{0}\|u\|_{L^{3,\infty}(B(\rho))}^{2}\|\Delta\phi+\partial_{t}\phi\|_{L_{3,1}(B(\rho))}ds\leqslant$$$$\leqslant c\frac{\rho}{(\rho-s)^{2}}\int\limits_{-\rho^{2}}^{0}\|u\|_{L^{3,\infty}(B(\rho))}^{2}ds\leqslant\frac{\rho^{\frac{5}{3}}}{(\rho-s)^{2}}(\int\limits_{-\rho^{2}}^{0}\|u\|_{L^{3,\infty}(B(\rho))}^{3}ds)^{\frac{2}{3}}.$$
Now we treat $(3)$.
Again by O'Neil  obtain:
\begin{equation}\label{est(3)}
(3)\leqslant 2\int\limits_{-\rho^2}^0\|u.\nabla\phi\|_{L^{3,1}(B(\rho))}\|p\|_{L_{\frac{3}{2},\infty}(B(\rho))}ds.
\end{equation}
By a well known  interpolation characterisation of Lorentz spaces:
$$L^{3,1}(B(\rho)=(L_{2}(B(\rho)),L_{6}(B(\rho))_{\frac{1}{2},1}.$$
Thus, by well known properties of interpolation spaces and Gagliardo Nirenberg inequality:
$$\|u.\nabla \phi\|_{L^{3,1}(B(\rho))}\leqslant\|u.\nabla \phi\|_{L_{2}(B(\rho))}^{\frac{1}{2}}\|u.\nabla \phi\|_{L_{6}(B(\rho))}^{\frac{1}{2}}\leqslant$$$$C\|u.\nabla \phi\|_{L_{2}(B(\rho))}^{\frac{1}{2}}\|\nabla(u.\nabla \phi)\|_{L_{2}(B(\rho))}^{\frac{1}{2}}\leqslant$$$$ \frac{C\|u\|_{L_{2}(B({\rho}))}}{(\rho-s)^{\frac{3}{2}}}+\frac{C\|u\|_{L_{2}(B({\rho}))}^{\frac{1}{2}}\|\nabla u\|_{L_{2}(B({\rho}))}^{\frac{1}{2}}}{\rho-s}.$$
Using this and applying Holder to (\ref{est(3)}) to get
$$(3)\leqslant \frac{C}{(\rho-s)^{\frac{3}{2}}}(\int\limits_{-\rho^{2}}^0 \|u\|_{L_{2}(B({\rho}))}^{3} ds)^{\frac{1}{3}}$$$$+\frac{C}{\rho-s}(\int\limits_{-\rho^{2}}^0 \|u\|_{L_{2}(B({\rho}))}^{\frac{3}{2}} \|\nabla u\|_{L_{2}(B({\rho}))}^{\frac{3}{2}}ds)^{\frac{1}{3}} \times (\int\limits_{-\rho^{2}}^{0}\|p\|_{L_{\frac{3}{2},\infty}(B(\rho))}^{\frac{3}{2}}ds)^{\frac{2}{3}}.$$
From here, it is not difficult to obtain
$$(3)\leqslant (\frac{Cr^{\frac{2}{3}}I_{1}(\rho)^{\frac{1}{2}}}{(\rho-s)^{\frac{3}{2}}}+\frac{Cr^{\frac{1}{6}}}{\rho-s}I_{2}(\rho)^{\frac{1}{4}}I_{1}(\rho)^{\frac{1}{4}})\times(\int\limits_{-r^{2}}^{0}\|p\|_{L^{\frac{3}{2},\infty}(B(r))}^{\frac{3}{2}}ds)^{\frac{2}{3}}.$$
Identical reasoning gives
$$(2)\leqslant (\frac{Cr^{\frac{2}{3}}I_{1}(\rho)^{\frac{1}{2}}}{(\rho-s)^{\frac{3}{2}}}+\frac{Cr^{\frac{1}{6}}}{\rho-s}I_{2}(\rho)^{\frac{1}{4}}I_{1}(\rho)^{\frac{1}{4}})\times(\int\limits_{-r^{2}}^{0}\|u\|_{L^{{3},\infty}(B(r))}^{{3}}ds)^{\frac{2}{3}}.$$
Thus, we see by Young's inequality that
\begin{equation}\label{iterativeineq}
I(s)\leqslant \frac{r^{\frac{5}{3}}}{(\rho-s)^{2}}(\int\limits_{-r^{2}}^{0}\|u\|_{L^{3,\infty}(B(r))}^{3}ds)^{\frac{2}{3}}+\frac{1}{2}I(\rho)+$$$$+(\frac{Cr^{\frac{4}{3}}}{(\rho-s)^{3}}+\frac{Cr^{\frac{1}{3}}}{(\rho-s)^{2}})\left[(\int\limits_{-r^{2}}^{0}\|u\|_{L^{{3},\infty}(B(r))}^{{3}}ds)^{\frac{4}{3}}+(\int\limits_{-r^{2}}^{0}\|p\|_{L^{\frac{3}{2},\infty}(B(r))}^{\frac{3}{2}}ds)^{\frac{4}{3}}\right].
\end{equation}
By Lemma \ref{iterativeGiusti} obtain
$$I({r}/{2})\leq r^{-\frac{1}{3}}(\int\limits_{-r^{2}}^{0}\|u\|_{L^{3,\infty}(B(r))}^{3}ds)^{\frac{2}{3}}+$$$$+Cr^{-\frac{5}{3}}(\int\limits_{-r^{2}}^{0}\|u\|_{L^{{3},\infty}(B(r)}^{{3}}ds)^{\frac{4}{3}}+(\int\limits_{-r^{2}}^{0}\|p\|_{L^{\frac{3}{2},\infty}(B(r))}^{\frac{3}{2}}ds)^{\frac{4}{3}}).$$
From here the conclusion is immediate. $\Box.$
\begin{remark}\label{energyflatboundary}
 After appropriate re labellings, the analogous estimate holds for $(u,p)$ suitable weak solution near the flat part of the boundary. 
\end{remark}
 Next, we define \begin{equation}\label{interiorpressure}
D_{1}(z_{0},r;p):=\frac{1}{r^{\frac{3}{2}}}\int\limits_{t_0-r^2}^{t_0}(\int\limits_{B(x_{0},r)}|\nabla p|^{\frac{9}{8}}dx)^{\frac{4}{3}} dt.
\end{equation}

\begin{pro}\label{Interiorversionlorentz}
Let $(u,p)$ be a suitable weak solution in $Q(z_{0},1)$ 
Then for $0<r<1$ the following  holds (c is some universal constant):
\begin{equation}\label{energyweaknorms1interior}
A(z_{0},{r}/{2};u)+ B(z_0,{r}/{2};u)\leqslant c( C_{\infty}(z_{0},r;u)^{\frac{4}{3}}+C^{+}_{\infty}(z_{0},r;u)^{\frac{2}{3}}+$$$$+D_{1}(z_{0},r;p)^{\frac{2}{3}}C_{\infty}(z_{0},r;u)^{\frac{1}{3}}).
\end{equation}
\end{pro}

\textbf{Proof} The set up is the same as Lemma \ref{energyboundLorentz}. The main difference is estimation of $$(3):=\int\limits_{Q(\rho)}2u.\nabla\phi p dxds.$$
Using the solenodial condition we can write:
$$(3):=\int\limits_{Q(\rho)}2u.\nabla\phi( p-[p]_{B(r)}) dxds.$$
We note the Poincare inequality
$$\|p(\cdot,t)-[p]_{B(r)}(t)\|_{L^{\frac{3}{2},1}(B(r))}\leqslant Cr^{\frac{1}{3}}\|\nabla p(\cdot,t)\|_{L_{\frac{9}{8}}(B(r))}.$$
Thus, by using O'Neil in space and Holder in time:
$$|(3)|\leqslant \frac{Cr^{\frac{1}{3}}}{\rho-s}(\int\limits_{-r^2}^0\|\nabla p\|_{L_{\frac{9}{8}}(B(r))}^{\frac{3}{2}}ds)^{\frac{2}{3}}
(\int\limits_{-r^2}^0\|u\|_{L_{3,\infty}(B(r))}^{3}ds)^{\frac{1}{3}}).$$ The remainder of the proof is the same as Lemma \ref{energyboundLorentz}. $\Box$\\

\textbf{Proof of Lemma \ref{interiorpressureLorentzdecay}}\\
Without loss of generality take $z_{0}=0$.
For a.a $t\in]-\rho^2,0[$, the pressure $p$ meets the equation in sense of distributions in $B(\rho)$:
$$\Delta p(\cdot,t)=-\rm{ div\,div}\,v(\cdot,t)\otimes v(\cdot,t).$$
Decompose the pressure so that
$$p=p_1+p_2,$$
where 
\begin{equation}\label{Rieszpres}
p_{1}(\cdot,t):=R_{i}R_{j}(\chi_{B(\rho)}v_{i}v_{j}(\cdot,t)).
\end{equation}
Here $R_{i}$ is Riesz operator and we adopt summation convention.
It is not difficult to notice that in $B(\rho)$:
\begin{equation}\label{p2harmonic}
\Delta p_{2}(\cdot,t)=0.
\end{equation}
For $p_1$ we have 
\begin{equation}\label{Rieszpresbound}
\|p_1(\cdot,t)\|_{L^{\frac{3}{2},\frac{q}{2}}(B(\rho))}^{\frac{3}{2}}\leqslant c\|v\otimes v(\cdot,t)\|_{L^{\frac{3}{2},\frac{q}{2}}(B(\rho))}^{\frac{3}{2}}\leqslant c\|v(\cdot,t)\|_{L^{3,q}(B(\rho))}^{3}.
\end{equation}
Let $0<r\leqslant\frac{\rho}{2}$.
Using O'Neil's inequality it is not difficult to see
$$ |[p_{1}]_{B(r)}|(t)\leqslant \frac{c\|p_{1}(\cdot,t)\|_{L^{\frac{3}{2},\frac{q}{2}}(B(r))}}{r^2}.$$
One can then easily obtain
$$D_{q}(0,r;p_1)\leq \frac{c}{r^2}(\int\limits_{-r^2}^{0}\|p_{1}(\cdot,t)\|_{L^{\frac{3}{2},\frac{q}{2}}(B(r))}^{\frac{3}{2}}dt)\leqslant$$$$\leqslant\frac{c}{r^2}(\int\limits_{-\rho^2}^{0}\|p_{1}(\cdot,t)\|_{L^{\frac{3}{2},\frac{q}{2}}(B(\rho))}^{\frac{3}{2}}dt)$$
Using (\ref{Rieszpresbound}), it is not difficult to see
\begin{equation}\label{decayestRieszp}
D_{q}(0,r,p)\leqslant c\left[\left(\frac{\rho}{r}\right)^{2}C_{q}(0,\rho,v)+D_{q}(0,r,p_{2})\right].
\end{equation}
Since $p_2$ is a harmonic function, we see that
$$\sup_{x\in B(r)}|p_2(x,t)-[p_2]_{B(r)}(t)|\leqslant cr\sup_{x\in B(\frac{\rho}{2})}|\nabla p_2(x,t)|\leqslant$$$$\leqslant
\frac{cr}{\rho^{4}}\int\limits_{B(\rho)}|p_{2}(x,t)-[p_2]_{B(\rho)}(t)|dx\leqslant$$$$\leqslant \frac{cr}{\rho^3}\|p_{2}(\cdot,t)-[p_2]_{B(\rho)}(t)\|_{L_{\frac{3}{2},\frac{q}{2}}(B({\rho}))}.$$
The last line follows from O'Neils inequality and fact that 
$$\|\chi_{\Omega}\|_{L_{p,q}(\Omega)}\leqslant C_{p,q}|\Omega|^{\frac{1}{p}}.$$
The remaining parts of the proof use this last fact, but are otherwise identical to that in \cite{S2005}.
The remaining details are an exercise for the reader. $\Box.$

Before proving the main statement we introduce some  notation:
\begin{equation}\label{uL3}
C(z_{0},r;u):= r^{-2}\int\limits_{Q(z_0,r)} |u|^3 dz,
\end{equation}
\begin{equation}\label{pressure3/2}
D(z_0,r;p):= r^{-2}\int\limits_{Q(z_0,r)} |p|^{\frac{3}{2}} dz.
\end{equation}
The following version of $\epsilon$- regularity criteria of Caffarelli-Kohn-Nirenberg will be important to us in the sequel:
\begin{lemma}\label{CKN}
Let $(u,p)$ be a suitable weak solution in $Q(R)$. Then there exists a universal constants $\epsilon_{0}$ and $c_{0k}$ (with $k=1,2\ldots$) with the following property.  Assume
\begin{equation}\label{smallnessconditionCKN}
C(0,R;u)+D(0,R;p)<\epsilon_{0}.
\end{equation}
then for any natural number $k$, $\nabla ^{k-1}u$ is Holder continuous in $\Bar{Q}(R/2)$ and the following bound is valid:
\begin{equation}\label{spatialsmoothnessCKN}
\max_{\Bar{Q}(R/2)}|\nabla^{k-1}u(z)|<c_{ok}R^{-k}.
\end{equation}

\end{lemma}
A short proof can be found in \cite{Lin} for example, a detailed one in \cite{Ladyser02}. Now we state our main result, here it is:

\begin{theorem}\label{strengthendepsilonregularityappendix}
Let $(u,p)$ be a suitable weak solution in $Q$. Then there exists a universal constants $\epsilon_{0}$ and $c_{0k}$ (with $k=1,2\ldots$) with the following property.  Assume
\begin{equation}\label{smallnessconditionappendix}
C_{\infty}(0,1;u)+D_{\infty}(0,1;p)<\epsilon_{0}.
\end{equation}
then for any natural number $k$, $\nabla ^{k-1}u$ is Holder continuous in $\Bar{Q}({1}/{8})$ and the following bound is valid:
\begin{equation}\label{spatialsmoothnessappendix}
\max_{\Bar{Q}({1}/{8})}|\nabla^{k-1}u(z)|<c_{ok}.
\end{equation}

\end{theorem}
\textbf{Proof}\\
From Lemma \ref{energyboundLorentz} and (\ref{smallnessconditionappendix}) it follows that
\begin{equation}\label{uniformenergy}
A(0,1/2;u)+B(0,1/2;u)\leqslant c\left(\epsilon_0+\epsilon_0^2\right)^{\frac 2 3}.
\end{equation}
Here $c$ will always denote some universal constant.
By interpolation and Sobolev embedding theorem one can show that
 $$
C(0,1/2,u)\leqslant c[A(0,1/2;u)^{\frac{3}{4}}B(0,1/2;u)^{\frac{3}{4}}+A(0,1/2;u)^{\frac{3}{2}}].
$$
Thus, by (\ref{uniformenergy}) we have
\begin{equation}\label{L3u}
C(0,1/2,u)\leqslant  c\left(\epsilon_0+\epsilon_0^2\right).
\end{equation}
For similar reasons it is not so difficult to see that 

$$
\|\rm{div}\,(u\otimes u)\|_{L_{\frac{9}{8},\frac{3}{2}}(Q(\frac{1}{2}))}\leqslant c[A(0,1/2;u)+A(0,1/2;u)^{\frac{1}{3}}B(0,1/2;u)^{\frac{2}{3}}].
$$
Thus,
\begin{equation}\label{nonlinear9/8,3/2}
\|\rm{div}\,(u\otimes u)\|_{L_{\frac{9}{8},\frac{3}{2}}(Q(\frac{1}{2}))}\leqslant \left(\epsilon_0+\epsilon_0^2\right)^{\frac 2 3}.
\end{equation}
By Holder's inequality it is obvious that
\begin{equation}\label{gradu9/8,3/2}
\|u\|_{W^{1,0}_{\frac 9 8,\frac 3 2}(Q(\frac{1}{2}))}\leqslant c(A(u,0;1/2)^{\frac 1 2}+B(u,0;1/2)^{\frac 1 2})\leqslant$$$$\leqslant c\left(\epsilon_0+\epsilon_0^2\right)^{\frac{1}{3}}
\end{equation}
Using O'Neil's inequality, we have
$$\int\limits_{B(1/2)} |p(x,t)|^{\frac 9 8} dx\leqslant c\||p(\cdot,t)|^{\frac 9 8}\|_{L^{\frac{8}{3},\infty}(B(\frac{1}{2}))}=
c\|p(\cdot,t)\|_{L^{3,\infty}(B(\frac{1}{2}))}^{\frac{9}{8}}.
$$
Hence, 
\begin{equation}\label{pres9/8,3/2}
\|p\|_{L_{\frac{9}{8},\frac{3}{2}}(Q(\frac{1}{2}))}\leqslant c \epsilon_{0}^{\frac{2}{3}}.
\end{equation}

Local interior regularity theory for Stokes equation gives

\begin{equation}\label{higherregrescailing}
\|\partial_{t}u\|_{L_{\frac{9}{8},\frac{3}{2}}(Q(\frac{1}{4}))}+\|\nabla^2u\|_{L_{\frac{9}{8},\frac{3}{2}}(Q(\frac{1}{4}))}+$$$$+\|\nabla p\|_{L_{\frac{9}{8},\frac{3}{2}}(Q(\frac{1}{4}))}\leqslant$$$$\leqslant
c[\|\rm{div}\,((u\otimes u)\|_{L_{\frac{9}{8},\frac{3}{2}}(Q(\frac{1}{2}))}+\|u\|_{L_{\frac{9}{8},\frac{3}{2}}(Q(\frac{1}{2}))}+$$$$+
\|\nabla u\|_{L_{\frac{9}{8},\frac{3}{2}}(Q(\frac{1}{2}))}+\|p\|_{L_{\frac{9}{8},\frac{3}{2}}(Q(\frac{1}{2}))}].
\end{equation} 
Using this together with (\ref{nonlinear9/8,3/2})- (\ref{pres9/8,3/2}) obtain that
$$\|\nabla p\|_{L_{\frac{9}{8},\frac{3}{2}}(Q(\frac{1}{4}))}\leqslant c[\left(\epsilon_0+\epsilon_0^2\right)^{\frac{1}{3}}+
\left(\epsilon_0+\epsilon_0^2\right)^{\frac{2}{3}}] .$$
Thus, by the Poincare inequality:
$$\|p-[p]_{,\frac{1}{4}}\|_{L_{\frac{3}{2}}(Q(\frac{1}{4}))}\leqslant c[\left(\epsilon_0+\epsilon_0^2\right)^{\frac{1}{3}}+
\left(\epsilon_0+\epsilon_0^2\right)^{\frac{2}{3}}] .$$
But by O'Neil $|[p(\cdot,t)]_{,\frac{1}{4}}|\leqslant c\|p(\cdot,t)\|_{L^{\frac{3}{2},\infty}(B)}$. Therefore, we conclude
\begin{equation}\label{pkuniform}
\|p\|_{L_{\frac{3}{2}} (Q(\frac{1}{4}))}\leqslant c[\left(\epsilon_0+\epsilon_0^2\right)^{\frac{1}{3}}+
\left(\epsilon_0+\epsilon_0^2\right)^{\frac{2}{3}}] .
\end{equation} 
This along with (\ref{L3u}) gives
$$ C(u,0;1/4)+D(u,0;1/4)\leqslant c[\left(\epsilon_0+\epsilon_0^2\right)^{\frac{1}{2}}+
\left(\epsilon_0+\epsilon_0^2\right)].$$
Choosing $\epsilon_{0}$ sufficiently small, gives the conclusion by Lemma \ref{CKN}.
 $\Box$
 \begin{remark}
  With certain adjustments to the proof, namely local boundary regularity for the Stokes system and the boundary analogue of \ref{CKN} (see  \cite{S3} and \cite{Ser09}), we have the following near the flat part of the boundary. \\\\
 Let $(u,p)$ be a suitable weak solution to the Navier-Stokes equations in $Q^+(1)$ near $\Gamma(0,1)\times ]-1,0[$. Then there exists universal constants $\epsilon_{0}$ and $c_0$ such that if
 $$C_{\infty}^{+}(0,1;u)+D_{\infty}^{+}(0,1;p)<\epsilon_{0},$$
 then $u$ is Holder continuous in $\bar{Q}^+(1/2)$ and the following is valid:
 $$\sup_{\bar{Q}^+(1/2)}|u(z)|\leqslant c_0.$$ 
 \end{remark}
\textbf{Acknowledgement.}\\
The Author is grateful to G. Seregin for helpful comments and
support.


\begin{thebibliography}{99}
\bibitem{Bergh}
Bergh, Jöran; Lofstrom, Jörgen Interpolation spaces. An introduction. Grundlehren der Mathematischen Wissenschaften, No. 223. Springer-Verlag, Berlin-New York, 1976.
\bibitem {CKN}
Caffarelli, L., Kohn, R.-V., Nirenberg, L., \emph{Partial regularity of
suitable weak solutions of the Navier-Stokes equations}, Comm. Pure
Appl. Math., Vol. XXXV (1982), pp. 771--831.

\bibitem{ESS2003}
Escauriaza, L.; Seregin, G.; \v Sver\'ak, V. $L_{3,\infty}$-solutions of Navier-Stokes equations and backward uniqueness. (Russian) Uspekhi Mat. Nauk 58 (2003), no. 2(350), 3--44; translation in Russian Math. Surveys 58 (2003), no. 2, 211–250.
\bibitem{Giusti}
 Giusti, Enrico Direct methods in the calculus of variations. World Scientific Publishing Co., Inc., River Edge, NJ, 2003. 
\bibitem{Hopf}
Hopf, Eberhard Über die Anfangswertaufgabe für die hydrodynamischen Grundgleichungen. (German) Math. Nachr.  4,  (1951). 213–231.

\bibitem{Serrin}
J. Serrin, On the interior regularity of weak solutions of the Navier–Stokes equations, Arch. Rational
Mech. Anal. 9 (1962), 187–195. MR0136885 (25:346)



\bibitem{kimkozono}
Kim, Hyunseok; Kozono, Hideo Interior regularity criteria in weak spaces for the Navier-Stokes equations. Manuscripta Math.  115  (2004),  no. 1, 85–100.
\bibitem{kiselevL} Kiselev, A. A.; Ladyženskaya, O. A. On the existence and uniqueness of the solution of the nonstationary problem for a viscous, incompressible fluid. (Russian) Izv. Akad. Nauk SSSR. Ser. Mat.  21  1957 655–680.
\bibitem{kozono}
 Kozono, Hideo Removable singularities of weak solutions to the Navier-Stokes equations. Comm. Partial Differential Equations  23  (1998),  no. 5-6, 949–966. 
\bibitem{L1967}
 Ladyženskaja, O. A. Uniqueness and smoothness of generalized solutions of Navier-Stokes equations. (Russian) Zap. Naučn. Sem. Leningrad. Otdel. Mat. Inst. Steklov. (LOMI)  5  1967 169–185. 
\bibitem{L1970}
Ladyzhenskaya, O. A., Mathematical problems of the dynamics of viscous 
incompressible fluids, 2nd edition, Nauka, Moscow 1970.
\bibitem{Ladyser02}
Ladyzhenskaya, O. A.; Seregin, G. A. On partial regularity of suitable weak solutions to the three-dimensional Navier-Stokes equations. J. Math. Fluid Mech. 1 (1999), no. 4, 356–387.
\bibitem {LSU}
Ladyzhenskaya, O. A., Solonnikov, V. A., Uralt'seva, N. N., Linear
and quasi-linear equations of parabolic type, Moscow, 1967;
English translation, American Math. Soc., Providence 1968.
\bibitem {Le}
J. Leray, \emph{ Sur le mouvement d'un liquide visqueux emplissant
l'espace}, Acta Math. \textbf{63} (1934),  193--248.

\bibitem{Lin}
Lin, Fanghua A new proof of the Caffarelli-Kohn-Nirenberg theorem. Comm. Pure Appl. Math. 51 (1998), no. 3, 241–257.

\bibitem{Tsai}
 Luo, Y., Tsai, T. P.,
Regularity criteria in weak $L_3$ for 3D incompressible Navier-Stokes equations.  arXiv:1310.8307, April 2014.

\bibitem{Struwe}
M. Struwe, On partial regularity results for the Navier–Stokes equations, Comm. Pure Appl. Math.
41 (1988), 437–458. MR0933230 (89h:35270)



 \bibitem{MSh2006}
 Mikhailov, A. S.; Shilkin, T. N. $L_{3,\infty}$  -solutions to the 3D-Navier-Stokes system in the domain with a curved boundary.  Zap. Nauchn. Sem.  (POMI)  336  (2006),  Kraev. Zadachi Mat. Fiz. i Smezh. Vopr. Teor. Funkts. 37, 133--152, 276;  translation in  J. Math. Sci. (N. Y.)  143  (2007),  no. 2, 2924–2935.
\bibitem{Phuc} Phuc,  N. C.. The Navier-Stokes equations in nonendpoint borderline Lorentz spaces. arXiv:1407.5129, July 2014.
\bibitem{Scheffer}
 Scheffer, V. Hausdorff measure and the Navier-Stokes equations. Comm. Math. Phys. 55 (1977), no. 2, 97–112. 
 
\bibitem{ser01}
 Seregin, G. A. On the number of singular points of weak solutions to the Navier-Stokes equations. Comm. Pure Appl. Math.  54  (2001),  no. 8, 1019–1028.
\bibitem{sersver02}
Seregin, G.; Sverak, V. Navier-Stokes equations with lower bounds on the pressure. Arch. Ration. Mech. Anal. 163 (2002), no. 1, 65–86.
\bibitem {S3}
Seregin, G.,\emph{ Local regularity of suitable weak solutions to the
Navier-Stokes equations near the boundary},  J. math. fluid mech.,
4(2002), no.1,1--29.


\bibitem{S2005}
Seregin, G. A. , \emph{ On smoothness of $L_{3,\infty}$-solutions to the Navier-Stokes
equations up to boundary}, Mathematische
Annalen, 332(2005), pp. 219-238.
\bibitem{Ser09}  Seregin, G. A.,   \emph{ A note on local boundary regularity for the Stokes system}, Zapiski Nauchn. Seminar., POMI, 370 (2009), pp. 151-159.
\bibitem{Ser12} Seregin, G. A., A certain necessary condition of potential blow up for Navier-Stokes equations. Comm. Math. Phys. 312 (2012), no. 3, 833–845. 

\bibitem{WangZhang1}
 Wang, W.; Zhang, Z. On the interior regularity criterion and the number of singular points to the Navier-Stokes equations.  arXiv:1201.1100, 5 Jan 2012. 
 \bibitem{WangZhang2}
 Wang, W.; Zhang, Z. Blow-up of critical norms for the 3-D Navier-Stokes equations.  arXiv:1510.02589 , 9 Oct 2015.
\end{thebibliography}
\end{document}